\documentclass[11pt]{amsart}

\usepackage{amscd}
\usepackage{amsmath}
\usepackage{graphicx}
\usepackage{amsfonts}
\usepackage{amssymb}
\begin{document}
\title[Characterizing derivations for any nest algebras on Banach spaces]
{Characterizing derivations for any nest algebras on Banach spaces
by their behaviors at an injective operator}
\author{Yanfang Zhang}
\address[Yanfang Zhang]{Department of Mathematics, Taiyuan University of Technology,
Taiyuan 030024, P. R. China.} \email{sxzyf1012@sina.com}

\author{Jinchuan Hou}
\address[Jinchuan Hou]{Department of Mathematics, Taiyuan University of Technology, Taiyuan 030024, P. R.  China}
\email{jinchuanhou@aliyun.com}

\author{Xiaofei Qi}
\address[Xiaofei Qi]{
Department of Mathematics, Shanxi University, Taiyuan 030006, P. R.
China.} \email{xiaofeiqisxu@aliyun.com}

\thanks{ This work is partially supported by National Natural Science Foundation of
China (11171249,11101250, 11271217) and Youth Foundation of Shanxi
Province (2012021004).}

\begin{abstract}

Let ${\mathcal N}$ be a nest on a  complex Banach space $X$ and let
$\mbox{\rm Alg}{\mathcal N}$ be the associated nest algebra. We say
that an operator $Z\in {\rm Alg}{\mathcal N}$ is an all-derivable
point of $\mbox{\rm Alg}{\mathcal N}$ if every linear map $\delta$
from $\mbox{\rm Alg}{\mathcal N}$ into itself derivable at $Z$ (i.e.
$\delta$ satisfies  $\delta(A)B+A\delta(B)=\delta(Z)$ for any $A,B
\in {\rm Alg}{\mathcal N}$ with $AB=Z$) is a derivation. In this
paper, it is shown that
 every injective operator and every operator
with dense range in ${\rm Alg}{\mathcal N}$ are all-derivable points
of ${\rm Alg}{\mathcal N}$ without any additional assumption on the
nest.

\end{abstract}
\maketitle

\section {Introduction}

Let $\mathcal A$ be an (operator) algebra. Recall that a linear map
$\delta:{\mathcal A}\rightarrow{\mathcal A}$  is  a derivation if
$\delta(AB)=\delta(A)B+A\delta(B)$ for all $A$, $B\in \mathcal A $.
As well known, the class of derivations is a very important class of
linear maps both in theory and applications, and was studied
intensively. The question of under what conditions that a linear
(even additive) map becomes a derivation attracted much attention of
authors (for instance, see \cite{Cr,K,LS,P} and the references
therein). One approach is to characterize  derivations by their
local behaviors.    We say that a map $\varphi:{\mathcal
A}\rightarrow{\mathcal A}$ is derivable at a point $Z\in {\mathcal
A}$ if $\varphi(A)B+A\varphi(B)=\varphi(Z)$ for any $A,B \in
{\mathcal A}$ with $AB=Z$, and we call such $Z$   a derivable point
of $\varphi$. Obviously,  a linear map is a derivation if and only
if it is derivable at every point. It is natural and interesting to
ask the question whether or not a linear map is a derivation if it
is derivable only at one given point. As usual, we say that an
element $Z\in {\mathcal A}$ is an all-derivable point of $\mathcal
A$ if every linear map on ${\mathcal A}$  derivable at $Z$ is in
fact a derivation. So far, we have known that  there exist many
all-derivable points (or full-drivable points) for certain
(operator) algebras (see \cite{HA,HQ3,LP,HQ1,Zh,ZZ} and the
references therein). However, zero point $0$ is not an all-derivable
point for any algebra because the generalized derivations are
derivable at $0$ \cite{HQ1}.

The unit $I$, and more generally, invertible elements, are
all-derivable points for many algebras and is a start point to find
other all-derivable points. For instance, $I$ is an all-derivable
point of prime rings and triangular algebras, and every invertible
element is an all-derivable point of $\mathcal J$-subspace lattice
algebras (see \cite{HA,HQ3,QH7}). As nest algebras are of an
important class of operator algebras, there are many papers on
finding all-derivable points of certain nest algebras. We mention
some results related to this paper. Zhu and Xiong \cite{JZ1} proved
that every strongly operator topology continuous linear map
derivable at $I$ between nest algebras on complex separable Hilbert
spaces is a
 derivation.
There they said that $I$ is an all-derivable point related to strong
operator topology. An and Hou in \cite{HA} generalized the above
result
and showed  that every linear map derivable at $I$ between nest
algebras on complex Banach spaces is a derivation, under the
additional assumption of the existence of a complemented nontrivial
element in the nest. In \cite{HQ1}, Qi and Hou showed further that,
if   the nest ${\mathcal N}$ on a Banach space satisfies ``$N\in
{\mathcal N}$ is complemented in the given Banach space whenever
$N_-=N$", then the unit operator $I$ is an all-derivable point of
the nest algebra Alg$\mathcal N$; furthermore, they show that every
injective operator and every operator with dense range in the nest
algebra are all-derivable points of the nest algebra. This
additional assumption on the nest is quite weak: such nest concludes
all nests on Hilbert spaces, all finite nest, all nest with
order-type $\omega+1$ or $1+\omega^{*}$ or $1+\omega^{*}+\omega+1$,
where $\omega$ is the order-type of natural numbers and $\omega^{*}$
is its anti-order-type. But the problem whether the unit operator
$I$, or every injective operator, or every operator with dense
range, is an all-derivable point of any nest algebras on any Banach
spaces remains open.

The purpose of the  present paper is to solve the above problem and
show that every injective operator and every operator with dense
range are all-derivable points for all nest algebras on  complex
Banach spaces without any additional assumptions on the nests.

The following is the main result of this paper.

{\bf Theorem 1.1.} {\it Let $\mathcal N$  be a   nest on a
 complex Banach space $X $ with $\dim X\geq 2$ and $\delta: \mbox{\rm Alg}{\mathcal
N}\rightarrow \mbox{\rm Alg}{\mathcal N}$ be a linear map.  Let
$Z\in{\rm Alg}{\mathcal N}$ be an injective operator or an
  operator with dense range in ${\rm Alg}{\mathcal N}$.  Then $\delta$ is derivable at the  operator $Z$ if and
only if $\delta$ is a derivation. That is, every injective operator
and every operator with dense range are all-derivable points of any
nest algebras.}

Particularly, we have

{\bf Corollary 1.2.} {\it Let $\mathcal N$  be a   nest on a complex
Banach space $X $ with $dim X\geq 2$. Then every invertible operator
in ${\rm Alg}{\mathcal N}$ is an all-derivable point of ${\rm
Alg}{\mathcal N}$.}

The paper is organized as follows. We fix some notations and
preliminary lemmas in Section 2 and give a proof of Theorem 1.1 in
Section 3.

\section {Preliminaries and lemmas}

In this section we fix some notations and give some lemmas.

Assume that $X$ is a Banach space over the complex field $\mathbb
C$.  Denote by $ {\mathcal B}(X)$ the algebra of all bounded linear
operators on $X$. The topological dual space of $X$ (i.e. the set of
all bounded linear functionals on $X$) is denoted by $X^*$. Let
$X^{**}$ be the second dual space of $X$. The map $\kappa:x\mapsto
x^{**}$, defined by $x^{**}(f)=f(x)$ for all $f$ in $X^*$, is  the
canonical embedding from $X$ into $X^{**}$. For any $T\in{\mathcal
B}(X)$, its Banach adjoint operator $T^*$  is the map from $X^*$
into $X^*$ defined by $(T^{*}f)(x)=f(Tx)$ for any $f\in X^*$ and
$x\in X$. If $f\in X^*$ and $x\in X$, the operator $x\otimes f$ on
$X$ is defined by  $(x\otimes f)(y)=f(y)x$ for any $y\in X$.
$x\otimes f$ is rank one whenever both $x$ and $f $ are nonzero, and
every rank one operator has this form. It is easily seen that
$(x\otimes f)^{*}=f\otimes x^{**}$ and $x^{**}T^{*}=(Tx)^{**}$. For
any non-empty subset $N\subseteq X$, denote by $N^\perp$ the
annihilator of $N$, that is, $N^{\perp}=\{f\in X^{*}:f(x)=0$ for
every $x\in N \}$. Some times we use $\langle x,f\rangle$ to present
the value $f(x)$ of $f$ at $x$. In addition, we use the symbols
${\rm ran}(T)$ and ${\rm ker}(T)$ for the range and the kernel of
operator $T$, respectively.

 Recall that a nest $\mathcal N$ in  $X$ is a chain
of closed (under norm topology) linear subspaces of $X$ containing
the trivial subspaces $\{0\}$ and $X $, which is closed under the
formation  of arbitrary closed linear span (denoted by $\bigvee$ )
and intersection (denoted by $ \bigwedge)$. $ {\mbox {\rm
Alg}}\mathcal N $ denotes the associated nest algebra, which is the
algebra of all operators $T$ in ${\mathcal B}(X)$ such that
$TN\subseteq N$ for every element $N\in \mathcal N $. When
${\mathcal N}\neq \{\{0\}, X\}$, we say that $\mathcal N $ is
nontrivial. It is clear that ${\mbox {\rm Alg}}\mathcal N=\mathcal
B(X)$ if $\mathcal N $ is trivial. For $N\in \mathcal N$, let
$N_{-}=\bigvee\{M\in \mathcal N \mid M \subset N\}$,
$N_{+}=\bigwedge \{M \in \mathcal N\mid N\subset M \}$ and $
N_{-}^{\perp}=(N_{-})^{\perp}$. Also, let $\{0\}_-=\{0\}$ and
$X_+=X$. Denote $\mathcal D_1(\mathcal N)=\bigcup \{N\in \mathcal
N\mid N_{-}\neq X\}$ and $\mathcal D_{2}(\mathcal N)=\bigcup
\{N_{-}^{\perp}\mid N\in \mathcal N \ \ {\mbox {\rm {and}} }\ \
N\neq \{0\}\}$. Note that $\mathcal D_{1}(\mathcal N)$ is dense in
$X$ and $\mathcal D_{2}(\mathcal N)$ is dense in $X^*$. Clearly,
$\mathcal D_{1}(\mathcal N)=X$ if and only if $X_-\not= X$ and
$\mathcal D_{2}(\mathcal N)=X^*$ if and only if $\{0\}\not=
\{0\}_+$.
  For more informations on nest algebras, we
refer to \cite{D,S}.

The following lemmas are needed to prove the main result.

{\bf Lemma 2.1.} (\cite{D,S}) {\it Let $\mathcal N$  be a nest on a
real or complex  Banach space $X $. The rank one operator $x\otimes
f$ belongs to ${\mbox {\rm Alg}}\mathcal N$ if and only if there is
some $N\in\mathcal N$ such that $x\in N$ and $f\in N_-^\perp$.}

{\bf Lemma 2.2.} {\it Let $\mathcal N$  be a   nest on a real or
complex Banach space $X $ with $\dim X\geq 2$. Assume that $\delta:
\mbox{\rm Alg}{\mathcal N}\rightarrow \mbox{\rm Alg}{\mathcal N}$ is
a linear map satisfying $\delta(P)=\delta(P)P+P\delta(P)$ for all
idempotent operators $P\in \mbox{\rm Alg}{\mathcal N}$ and
$\delta(N)N+N\delta(N)=0$ for all operators $N\in \mbox{\rm
Alg}{\mathcal N}$ with $N^2=0$. If $ X_-\not= X$, then,}

(1) {\it for any $x\in X$ and $f\in X_{-}^{\perp}$, we have
$\delta(x\otimes f)\ker (f)\subseteq \mbox{\rm span}\{x\}$;}

(2) {\it there exist  linear transformations $B:X\rightarrow X$ and
 $C:X_{-}^{\perp}\rightarrow X^{\ast}$ such that $\delta(x\otimes
f)=Bx\otimes f+x\otimes Cf $ and $\langle Bx,f \rangle +\langle
x,Cf\rangle=0 $ for all $x\in X$ and $f\in X_{-}^{\perp}$.}

{\bf Proof.} (1) For any nonzero $x\in X$ and $f\in X_{-}^{\perp}$,
it follows from Lemma 2.1 that
 $x\otimes f\in {\rm Alg} \mathcal {N}$. If $\langle x,f \rangle\neq 0$,
letting $\bar{x}=\langle x, f \rangle^{-1} x$, then $\bar{x}\otimes
f $ is an idempotent operator in ${\rm Alg}\mathcal{N}$. By the
assumption on $\delta$, we have
$$\delta(\bar{x}\otimes f)=\delta(\bar{x}\otimes f)(\bar{x}\otimes
f)+(\bar{x}\otimes f)\delta(\bar{x}\otimes f).\eqno (2.1)$$ For any
$y\in \ker (f)$, letting the operators in Eq.(2.1) act at $y$, we
obtain  $\delta(\bar{x}\otimes f)y=\langle \delta(\bar{x}\otimes f)y
,f \rangle \bar{x}$, which implies $\delta(x\otimes
f)y=\langle\delta(\bar{x}\otimes f)y,f\rangle x\in {\rm span}\{x\}.$
If $ \langle x,f \rangle=0$, we can take  $z\in X$ such that
$\langle z,f\rangle=1$ as $X_-\not=X$. It is obvious  that both
$(x+z)\otimes f$ and $z\otimes f $ are idempotents   in ${\rm
Alg}\mathcal {N}$. By the assumption on $\delta$ again,  one has
$$\begin{array}{rl}&\delta(x\otimes f)=\delta((x+z)\otimes
f)-\delta(z\otimes f) \\
=&\delta(x\otimes f)(z\otimes f)+(x\otimes f)\delta(x\otimes
f)+\delta (z\otimes f)(x\otimes f)+(z\otimes f)\delta(x\otimes
f)\\
\end{array}\eqno(2.2)$$ and $$0=\delta(x\otimes f)(x\otimes f)+(x\otimes f)\delta (x\otimes
f).\eqno(2.3)$$ Then, for any $y\in \ker (f)$, Eqs.(2.2)-(2.3)
become to $\delta(x\otimes f)y=\langle\delta(z\otimes f)y,f\rangle
x+\langle\delta(x\otimes f)y,f \rangle z$ and
$\langle\delta(x\otimes f)y,f \rangle=0$. So $\delta(x\otimes
f)y=\langle\delta(z\otimes f)y,f\rangle x$ holds for every $y\in
\ker(f)$, that is, (1) holds.

(2) For any $x\in X$ and $f\in X_{-}^{\perp}$, by (1), there exists
a functional $\varphi_{x,f}$ on $\ker (f)$ such that
$\delta(x\otimes f)y=\varphi_{x,f}(y)x$ for all $y\in \ker(f)$. It
is easy to see that $\varphi_{x,f}$   is linear. Take a non-zero
vector $ w_f \in X\setminus X_{-} $ such that $ \langle w_f,f
\rangle =1 $. Let $\bar{\varphi}_{x,f}$ be a   linear extension of
$\varphi_{x,f}$ to $X$. Then
$\tilde{\varphi}_{x,f}=\bar{\varphi}_{x,f}-\bar{\varphi}_{x,f}(w_f)f$
is also a   linear extension of $\varphi_{x,f}$ which vanishes at
$w_f$. Define a map $B_f:X\rightarrow X$ by $B_{f}x=\delta(x\otimes
f)w_f$. Obviously, $B_{f}$ is linear by the linearity of $\delta$,
and $B_{\lambda f}= B_f$ for any nonzero scalar $\lambda$. For any
$\tilde{x}\in X$, as $z_{x}=\tilde{x}-f(\tilde{x})w_f\in \ker(f)$,
we have
$$\begin{array}{rl}\delta(x\otimes f)\tilde x=&\delta(x\otimes
f)(z_{x}+f(\tilde{x})w_f)
=\delta(x\otimes f)z_{x}+f(\tilde{x})\delta(x\otimes f)w_f\\
=&\varphi_{x,f}(z_{x})x+f(\tilde{x})\delta(x\otimes f)w_f
=\varphi_{x,f}(\tilde{x}-f(\tilde{x})w_f)+f(\tilde{x})B_{f}x\\
=&\tilde{\varphi}_{x,f}(\tilde{x})x+f(\tilde{x})B_{f}x. \end{array}
$$ So
$$\delta(x\otimes f)=x\otimes \tilde{\varphi}_{x,f}+B_{f}x\otimes f
\quad{\rm for \ \ all}\quad x\in X,\ f\in X_-^{\perp}. \eqno(2.4) $$
As $\delta(x\otimes f)$ and $B_{f}x\otimes f$ are bounded linear
operators on $X$, we see that $\tilde{\varphi}_{x,f}\in X^*$.

For simplicity seek, in the following we still denote
${\varphi}_{x,f}$ for $\tilde{\varphi}_{x,f}$.

{\bf Claim 1.} $\varphi_{x,f} $ only depends on $f$.

Take any $x_{1},x_{2}\in X$ with $x_1\not=0$. For any $y\in X$, by
Eq.(2.4), we have
$$\delta((x_{1}+x_{2})\otimes f)y=\langle
y,\varphi_{x_{1}+x_{2},f}\rangle(x_{1}+x_{2})$$and
$$\delta((x_{1}+x_{2})\otimes f)y=\delta(x_{1}\otimes f)y+\delta(x_{1}\otimes f)y
=\langle y,\varphi_{x_{1},f}\rangle x_{1} +\langle
y,\varphi_{x_{2},f}\rangle x_{2}.$$ So$$0=(\langle
y,\varphi_{x_{1}+x_{2},f}\rangle-\langle
y,\varphi_{x_{1},f}\rangle)x_{1}+(\langle
y,\varphi_{x_{1}+x_{2},f}\rangle-\langle
y,\varphi_{x_{2},f}\rangle)x_{2}.$$   If $x_{1}$ and $x_{2}$ are
linearly independent, the above equation implies
$\varphi_{x_{1},f}=\varphi_{x_{1}+x_{2},f}=\varphi_{x_{2},f}$. If
$x_{1}$ and $x_{2}$ are linearly dependent, we take $x_{3}\in X$
such that it is independent of $x_{1}$. By the preceding proof, one
gets $\varphi_{x_{1},f}=\varphi_{x_{3},f}=\varphi_{x_{2},f}$. The
claim is true.

Now we can write $\varphi_{x,f}=\varphi_{f}\in X^*$. Then
$$\delta(x\otimes f)=x\otimes \varphi_{f}+B_{f}x\otimes f\quad{\rm for \ \ all}\quad x\in X,\ f\in X_-^{\perp}.\eqno(2.5)$$

{\bf Claim 2.} There exist  linear transformations $B:X\rightarrow
X$ and   $C:X_{-}^{\perp}\rightarrow X^{\ast}$ such that
$\delta(x\otimes f)=Bx\otimes f+x\otimes Cf $ for all $x\in X$ and
$f\in X_{-}^{\perp}$.

Fix a non-zero $f_{0}\in X_{-}^{\perp}$ and put $B=B_{f_{0}}.$

For any nonzero $f_{1},\ f_{2}\in X_{-}^{\perp}$,  we show  that the
difference $B_{f_{1}}-B_{f_{2}}$ is a scalar multiple of the
identity $I$. If $f_{1}$ and $f_{2}$ are linearly independent, then
there are $x_{1},x_{2}\in X$ such that $ f_{i}(x_{i})=1$ and
$f_{i}(x_{j})=0$ for $1\leq i\neq j\leq 2$. For any $x\in X$, we
have
$$\delta(x\otimes(f_{1}+f_{2}))=x\otimes
\varphi_{f_{1}+f_{2}}+B_{f_{1}+f_{2}}x\otimes(f_{1}+f_{2}) $$ and
$$\begin{array}{rl}\delta(x\otimes(f_{1}+f_{2}))=&\delta(x\otimes
f_{1})+\delta(x\otimes f_{2})\\
=&x\otimes (\varphi_{f_{1}}+ \varphi_{f_{2}})+B_{f_{1}}x\otimes
f_{1}+B_{f_{2}}x\otimes f_{2}.\end{array}$$ So
$$0=x\otimes (\varphi_{f_{1}+f_{2}}-\varphi_{f_{1}}-\varphi_{f_{1}})+B_{f_{1}}x\otimes f_{2}+B_{f_{2}}x\otimes f_{1}.\eqno(2.6)$$
By acting at  $x_{1}-x_2$,  Eq.(2.6) gives
$$B_{f_{1}}x-B_{f_{2}}x=(\varphi_{f_{1}+f_{2}}(x_{1}-x_{2})-\varphi_{f_{1}}(x_{1}-x_{2})-\varphi_{f_{1}}(x_{1}-x_{2}))x=\lambda_{f_1,f_2}x.$$
Thus,  $B_{f_{1}}-B_{f_{2}}=\lambda_{f_1,f_2}I$ is a scalar multiple
of the identity $I$.  If $f_{1}$ and $f_{2}$ are linearly dependent,
then $B_{f_{1}}-B_{f_{2}}=0$ as $B_{f_{1}}=B_{f_{2}}$. By now we
have shown that, for any $f\in X_{-}^{\perp}$,
$B_{f}=B_{f_0}+b_{f}I=B+b_fI$ for some scalar $b_{f}$.
 Then Eq.(2.5) becomes
$$\begin{array}{rl}\delta(x\otimes f )=&x\otimes \varphi_{f}+(B+b_{f}I)x\otimes f\\
=&Bx\otimes f+x\otimes(b_{f}f+\varphi_{f})=Bx\otimes f+x\otimes
Cf,\end{array}$$ where $C: X_-^{\perp}\rightarrow X^{\ast}$ is
defined by $Cf=b_{f}f+\varphi_{f}$ for all $f\in X_-^\perp$. By the
linearity of $\delta$, we see that $C$ is linear.

{\bf Claim 3.} $\langle Bx,f \rangle +\langle x,Cf\rangle=0 $ for
all $x\in X$ and $f\in X_{-}^{\perp}$.

For any $x\in X$ and $f\in X_{-}^{\perp}$, by the assumptions on
$\delta$, one can get $(x\otimes f)\delta(x\otimes f)(x\otimes
f)=0$. Then, by Claim 2,  $(\langle Bx,f \rangle +\langle
x,Cf\rangle)(x\otimes f)=0$, and so $\langle Bx,f \rangle +\langle
x,Cf\rangle=0 $.

Combining Claims 1-3, (2) is true. \hfill$\Box$

{\bf Lemma 2.3.} {\it Let $\mathcal N$  be a   nest on a real or
complex Banach space $X $ with $\dim X\geq 2$. Assume that $\delta:
\mbox{\rm Alg}{\mathcal N}\rightarrow \mbox{\rm Alg}{\mathcal N}$ is
a linear map satisfying $\delta(P)=\delta(P)P+P\delta(P)$ for all
idempotent operators $P\in \mbox{\rm Alg}{\mathcal N}$ and
$\delta(N)N+N\delta(N)=0$ for all operators $N\in \mbox{\rm
Alg}{\mathcal N}$ with $N^2=0$. If $\{0\}\not=\{0\}_{+}$, then}

(1) {\it for any $f\in X^{\ast}$ and any $x \in \{0\}_+$, we have
$\delta(x\otimes f)^{\ast}(\ker x^{\ast\ast})\subseteq {\rm
span}\{f\}$;}

(2) {\it there exist linear transformations $B:X^*\rightarrow X^*$
and   $C:\kappa(\{0\}_{+})\rightarrow X^{\ast\ast}$ such that
$\delta(x\otimes f)^{\ast}=Bf\otimes x^{\ast\ast}+f\otimes
Cx^{\ast\ast}$ and $\langle Bf,x^{\ast\ast}\rangle+\langle
f,Cx^{\ast\ast} \rangle=0$ holds for all $x\in\{0\}_+$ and $f\in
X^*$.}

{\bf Proof.}  This is the ``dual''   of Lemma 2.2. We   give a proof
of the conclusion (1) in detail as a sample, and give a proof of the
conclusion (2) in  sketch.

(1) Let $f\in X^{\ast}$  and $x \in \{0\}_+$ be nonzero. If $\langle
f,x^{\ast\ast}\rangle\neq 0$, it is enough to consider the case
 $\langle f,x^{\ast\ast}\rangle= \langle
x,f\rangle=1 $. In this case, $(x\otimes f)^2=x\otimes f$. By the
assumptions on $\delta$, we get
$$\delta(x\otimes
f)^{\ast}=(f\otimes x^{\ast\ast})\delta(x\otimes
f)^{\ast}+\delta(x\otimes f)^{\ast}(f\otimes
x^{\ast\ast}).\eqno(2.7)$$ Applying Eq.(2.7) to any $g\in
\ker(x^{\ast\ast})$  gives $\delta(x\otimes
f)^{\ast}g=\langle\delta(x\otimes f)^{\ast}g,x^{\ast\ast}\rangle
f\in{\rm span}\{f\} $.

If $\langle f,x^{\ast\ast}\rangle=\langle x,f\rangle=0 $,  we can
find $f_{1}\in X^*$ with $\langle x,f_{1}\rangle=1$ as $\{0\}\not=
\{0\}_{+}$. Then both $x\otimes f_{1}$ and $x\otimes (f_{1}+f) $ are
idempotents in $\mbox{\rm Alg}{\mathcal N}$ and $(x\otimes f)^2=0$.
So
$$\delta(x\otimes(f+f_{1}))=\delta(x\otimes(f+f_{1}))(x\otimes(f+f_{1}))+(x\otimes(f+f_{1}))\delta(x\otimes(f+f_{1})),$$
$$0=\delta(x\otimes f)(x\otimes f)+(x\otimes f)\delta(x\otimes
f)\eqno(2.8)$$and$$\delta(x\otimes f_{1})=\delta(x\otimes
f_{1})(x\otimes f_{1})+(x\otimes f_{1})\delta(x\otimes f_{1}).$$ The
above two equations imply
$$\begin{array}{rl}\delta(x\otimes f)^{\ast}=&(f_{1}\otimes
x^{\ast\ast})\delta(x\otimes f)^{\ast}+(f\otimes
x^{\ast\ast})\delta(x\otimes f_{1})^{\ast}\\&+\delta(x\otimes
f_{1})^{\ast}(f\otimes x^{\ast\ast})+\delta(x\otimes
f)^{\ast}(f_{1}\otimes x^{\ast\ast}).\end{array}\eqno( 2.9)$$ Note
that Eq.(2.8) implies $0=(f\otimes x^{\ast\ast})\delta(x\otimes
f)^{\ast}+\delta(x\otimes f)^{\ast}(f\otimes x^{\ast\ast}).$ For any
$g\in \ker(x^{\ast\ast})$, letting the equation act at   $g$, one
gets $\langle\delta(x\otimes f)^{\ast}g,x^{\ast\ast}\rangle f=0$,
and so $\langle\delta(x\otimes f)^{\ast}g,x^{\ast\ast}\rangle=0$.
Thus,  letting Eq.(2.9) act  at any $g\in \ker(x^{\ast\ast})$ leads
to $\delta(x\otimes f)^{\ast}g=\langle\delta(x\otimes
f_{1})^{\ast}g,x^{\ast\ast}\rangle
 f\in{\rm span}\{f\}$.

(2) For any nonzero $f\in X^*$ and $x\in \{0\}_+$, we have $x\otimes
f\in\mbox{\rm Alg}{\mathcal N}$. Take $h_x \in X^{\ast}$ such that
$\langle x,h \rangle=1$. Define a map $B_{x}:X^{\ast}\rightarrow
X^{\ast}$ by $B_{x}f=\delta(x\otimes f)^{\ast}h_x$. Let $\kappa:
X\rightarrow X^{\ast\ast}$ be the canonical map from $X$ into
$X^{\ast\ast}$. By (1) in the lemma, there exists a linear
functional $\Phi_{f,x}$ on $ \ker x^{\ast\ast}$ such that
$\delta(x\otimes f)^{\ast} g=\Phi_{x,f}(g)f$ for any $ g\in
\ker(x^{\ast\ast})$. Still by $\Phi_{x,f}$ denotes such a special
extension of $\Phi_{x,f}$ on $X^{\ast\ast}$ which vanishes at $h_x$.
For any $\tilde{f}\in X^{\ast}$, let $Z_{f}=\tilde{f}-\langle
f,x^{\ast\ast}\rangle h_x$. It is obvious that $Z_{f}\in \ker
x^{\ast\ast}$, and so
$$\delta(x\otimes f)^{\ast}\tilde{f}=\delta(x\otimes f)^{\ast}(Z_{f}+\langle f,x^{\ast\ast}\rangle h_x)
=\Phi_{x,f}(\tilde{f})f+(B_{x}f\otimes x^{\ast\ast})\tilde{f}.$$
Thus we can get $\delta(x\otimes f)^{\ast}=B_{x}f\otimes
x^{\ast\ast}+f\otimes \Phi_{x,f}$ for $f\in X^*$ and $x\in \{0\}_+$.

Similarly to the proof of Lemma 2.2(2), the following Claims 1-3
hold.

{\bf Claim 1.} $\Phi_{x,f}$ depends only on $x$.

 Fix a non-zero vector $x_{0}\in \{0\}_{+}$ and put
$B=B_{x_{0}}$.

{\bf Claim 2.} There exist   linear transformations
$B:X^*\rightarrow X^*$ and  $C:\kappa(\{0\}_{+})\rightarrow
X^{\ast\ast}$ such that $\delta(x\otimes f)^{\ast}=Bf\otimes
x^{\ast\ast}+f\otimes Cx^{\ast\ast}$ holds for all $x\in\{0\}_+$ and
$f\in X^*$.

{\bf Claim 3.} $\langle Bf,x^{\ast\ast}\rangle+\langle
f,Cx^{\ast\ast} \rangle=0$ holds for all $x\in \{0\}_{+}$ and $f\in
X^*$.

Hence (2) holds, the proof of the lemma is finished. \hfill$\Box$

{\bf Lemma 2.4.} {\it Let $X$ be  a Banach space over the real or
complex field $\mathbb F$. Suppose that $\mathcal N$ is a   nest on
$X $ and $\delta: \mbox{\rm Alg}{\mathcal N}\rightarrow \mbox{\rm
Alg}{\mathcal N}$ is a linear map satisfying
$\delta(P)=\delta(P)P+P\delta(P)$ for all idempotent operators $P\in
\mbox{\rm Alg}{\mathcal N}$ and $\delta(N)N+N\delta(N)=0$ for all
operators $N\in \mbox{\rm Alg}{\mathcal N}$ with $N^2=0$. If
$X_{-}=X$ and $\{0\}=\{0\}_{+}$, then }

(1) {\it there exists a bilinear functional $\beta: ({\mathcal
D}_1({\mathcal N})\times {\mathcal D}_2({\mathcal N}))\cap{\rm
Alg}{\mathcal N}\rightarrow \mathbb F$ such that $(\delta(x\otimes
f)-\beta_{x,f}I)\ker (f)\subseteq \mbox{\rm span}\{x\}$ holds for
all $x\otimes f\in {\rm Alg} \mathcal {N}$;}

(2) {\it  there exist linear transformations $B: \mathcal
D_1(\mathcal N)\rightarrow \mathcal D_1(\mathcal N)$ and $C:
\mathcal D_{2}(\mathcal N)\rightarrow \mathcal D_{2}(\mathcal N)$
such that $\delta(x\otimes f)-\beta_{x,f}I=x\otimes Cf+Bx\otimes f$
holds for all $x\otimes f\in {\mbox {\rm Alg}}\mathcal N$.}

{\bf Proof.} Since $X_{-}=X$ and $\{0\}=\{0\}_{+}$, it is obvious
that $\mathcal N$ is non-trivial, $\mathcal D_{1}(\mathcal N)$ and
$\mathcal D_{2}(\mathcal N)$ are dense proper linear manifolds in
$X$ and $X^*$, respectively, and, for each nontrivial $N\in{\mathcal
N}$, both $N$ and $N_-^\perp$ are infinite dimensional.

For any $x\otimes f\in {\rm Alg} \mathcal {N}$, if $\langle
x,f\rangle\neq 0$,  by the assumption on $\delta$ for idempotents,
it is easily seen that $\delta(x\otimes f)\ker(f)\in {\rm
span}\{x\}$. In this case let $\beta_{x,f}=0$.

Now assume that $x\otimes f\in$Alg$\mathcal N$ and $\langle
x,f\rangle = 0$. By Lemma 2.1, there exists some $N_{x}\in \mathcal
N$ such that $x\in N_{x}$ and $f\in (N_{x})_{-}^{\perp}$. It is easy
to check that $\vee\{\ker(f)\cap N:N\in \mathcal N, N_x\subseteq N,
N_-\not= X\}={\mathcal D}_1({\mathcal N})\cap\ker(f)$ is dense in
$\ker(f)$. For any $N\in\mathcal N$ with $N_x\subseteq N$ and for
every   $y\in N\cap \ker(f)$, take $g\in N^{\perp}$ such that $g$ is
linearly independent of $f$. Then $x\otimes f+y\otimes g\in
\mbox{\rm Alg}{\mathcal N}$ with $(x\otimes f+y\otimes g)^2=0$. By
the assumption on $\delta$, we have
$$0=\delta(x\otimes f+y\otimes g)(x\otimes f+y\otimes g)+(x\otimes f+y\otimes g)\delta(x\otimes f+y\otimes g).$$
Also note that $(x\otimes f)^2=0$ and $(y\otimes g)^2=0$. The above
equation can be reduced to
$$0=\delta(x\otimes f)(y\otimes g)+\delta(y\otimes g)(x\otimes f)+(y\otimes g)\delta(x\otimes f)+(x\otimes f)\delta(y\otimes g).\eqno( 2.10)$$
Choose $z_1\in X$ such that $\langle z_1,f\rangle=-1$ and  $\langle
z_1,g\rangle=1$. Let Eq.(2.10) act at $z_1$, one gets
$$0=\delta(x\otimes f)y-\delta(y\otimes g)x+\langle \delta(x\otimes f)z_1,g\rangle y+\langle \delta(y\otimes g)z_1,f\rangle x. \eqno( 2.11)$$
Since $f$ and $g$ are linearly independent, we can pick $w\in X$
such that $\langle w,f\rangle=2$ and $\langle w,g\rangle=0$. Let
$z_2=z_1+w$; then $\langle z_2,g\rangle=\langle z_2,f\rangle=1$.
Letting each operator in Eq.(2.10) act at $z_2$, we get
$$0=\delta(x\otimes f)y+\delta(y\otimes g)x+\langle \delta(x\otimes
f)z_2,g\rangle y+\langle \delta(y\otimes g)z_2,f\rangle x.
\eqno(2.12)$$ It follows from Eqs.(2.11)-(2.12) that
$$\delta(x\otimes f)y=-\frac{1}{2}(\langle \delta(x\otimes f)z_1,g\rangle+\langle \delta(x\otimes f)z_2,g\rangle) y
-\frac{1}{2}(\langle \delta(y\otimes g)z_1,f\rangle+\langle
\delta(y\otimes g)z_2,f\rangle )x,$$ which implies that
$\delta(x\otimes f)y\in {\rm span}\{x,y\},$ that is, there exist
scalars $\alpha_{x,f}^{N}(y)$, $\beta_{x,f}^{N}(y)$ such that
$$\delta(x \otimes
f)y=\alpha_{x,f}^{N}(y)x+\beta_{x,f}^{N}(y)y.
 \eqno(2.13)$$

Next we show that, for any $y_1,y_2\in N\cap\ker(f)$,
$\beta_{x,f}^{N}(y_1)=\beta_{x,f}^{N}(y_2)$. In fact, by Eq.(2.13),
we have
$$\delta(x \otimes
f)y_{i}=\alpha_{x,f}^{N}(y_{i})x+\beta_{x,f}^{N}(y_{i})y_{i},\quad\quad
i=1,2$$ and $$\delta(x \otimes
f)(y_{1}+y_{2})=\alpha_{x,f}^{N}(y_{1}+y_{2})x+\beta_{x,f}^{N}(y_{1}+y_{2})(y_{1}+y_{2}).$$
So
$$\begin{array}{rl}0=&(\alpha_{x,f}^{N}(y_{1}+y_{2})-\alpha_{x,f}^{N}(y_{1})-\alpha_{x,f}^{N}(y_{2}))x
\\
&+(\beta_{x,f}^{N}(y_{1}+y_{2})-\beta_{x,f}^{N}(y_{1}))y_{1}
+(\beta_{x,f}^{N}(y_{1}+y_{2})-\beta_{x,f}^{N}(y_{2}))y_{2}.\\\end{array}$$
If $y_1$, $y_2$ and $x$ are linearly independent, obviously
$\beta_{x,f}^{N}(y_{1})=\beta_{x,f}^{N}(y_{2});$ if $y_1$ and $y_2$
are linearly dependent and $x\not\in{\rm span}\{y_1,y_2\}$, as ${\rm
dim}(\ker f\cap N)=\infty$, there exists $y_3\in \ker f\cap N$ so
that $y_3\not\in{\rm span}\{y_1,y_2\}$. Then, by what  just proved
above, we have
$\beta_{x,f}^{N}(y_{1})=\beta_{x,f}^{N}(y_{3})=\beta_{x,f}^{N}(y_{2}).$
For the case $\dim{\rm span}\{x,y_1,y_2\}=1$, choosing  $y_4\in \ker
(f)\cap N$ so that $y_4$ is linearly independent of  $x$, we see
that
$\beta_{x,f}^{N}(y_{1})=\beta_{x,f}^{N}(y_{2})=\beta_{x,f}^{N}(y_{4})$.
Thus we have shown that $\beta_{x,f}^{N}(y)$ is independent of $y$.
So, there exists a scalar $\beta_{x,f}^{N}$ such that
$\beta_{x,f}^{N}(y)=\beta_{x,f}^{N}$ holds for all $y\in
N\cap\ker(f)$.

We  claim that $\alpha_{x,f}^{N}(y)$ and $\beta_{x,f}^{N}$ are
independent of $N$, and thus $\alpha_{x,f}^{N}(y)=\alpha_{x,f}(y)$,
$\beta_{x,f}^{N}=\beta_{x,f}.$ Indeed, if  ${N}^{'},{N}^{''}\in
\mathcal N$ with $N_x\subseteq N^{'}\cap N^{\prime\prime}$ and
$y\in{N}^{'}\cap N^{''}\cap\ker(f)$, then by Eq.(2.13), we have
$$\delta(x \otimes
f)y=\alpha_{x,f}^{{N}^{'}}(y)x+\beta_{x,f}^{{N}^{'}}y\quad{\rm
and}\quad\delta(x \otimes
f)y=\alpha_{x,f}^{{N}^{''}}(y)x+\beta_{x,f}^{{N}^{''}}y.$$ The above
two equations give
$$0=(\alpha_{x,f}^{{N}^{'}}(y)-\alpha_{x,f}^{{N}^{''}}(y))x+(\beta_{x,f}^{{N}^{'}}-\beta_{x,f}^{{N}^{''}})y.$$
As $\dim N^\prime\cap N^{\prime\prime}\cap\ker(f)=\infty$, one can
choose $y$ such that $y$ is linearly independent of $x$. Then we get
$\alpha_{x,f}^{{N}^{'}}(y)=\alpha_{x,f}^{{N}^{''}}(y)$ and
$\beta_{x,f}^{{N}^{'}}=\beta_{x,f}^{{N}^{''}}$.  It follows that
$\beta_{x,f}^{{N}}$ is independent of $N$ and   there is a scalar
$\beta_{x,f}$ such that $\beta_{x,f}^{{N}}=\beta_{x,f}$ holds for
all $N\in{\mathcal N}$ with $N_x\subseteq N$. Now it is clear that
$\alpha_{x,f}^{{N }^{'}}(y)=\alpha_{x,f}^{{N }^{''}}(y)$ also holds
for all
 $y\in N^\prime\cap N^{\prime\prime}\cap\ker(f)$. Hence there exists
a scalar $\alpha_{x,f}(y)$ such that
 $ \alpha_{x,f}^{{N}}(y)=\alpha_{x,f}(y)$ for all $y\in N\cap\ker(f)$ with $N_x\subseteq N$.

Thus we have shown that
  $$\delta(x \otimes
f)y=\alpha_{x,f}(y)x+\beta_{x,f}y
 \eqno(2.14)$$  holds  for all $ y\in
{\mathcal D}_1({\mathcal N})\cap\ker(f)$. Since ${\mathcal
D}_1({\mathcal N})\cap\ker(f)$ is dense in
 $\ker(f)$, Eq.(2.14) holds for all $y\in\ker (f)$.

Finally we will prove that $\beta_{x,f}$ is bilinear in
$x\in{\mathcal D}_1({\mathcal N})$ and $f\in{\mathcal D}_2({\mathcal
N})$ with $x\otimes f\in$Alg$\mathcal N$. Indeed, for any
$x_1,x_2\in{\mathcal D}_1({\mathcal N})$, there exists some
$N\in{\mathcal N}$ such that $x_{1},\ x_{2}\in
 N$. Then, for any $f\in N_{-}^{\perp}$, $x_{1}\otimes f$, $x_{2}\otimes f\in {\rm Alg}\mathcal {N}$.
By Eq.(2.14), for any $y\in\ker (f)$, one has
$$\delta((x_{1}+x_{2})\otimes
f)y=\alpha_{x_{1}+x_{2},f}(y)(x_{1}+x_{2})+\beta_{x_{1}+x_{2},f}y,$$
$$\delta(x_{1}\otimes
f)y=\alpha_{x_{1},f}(y)(x_{1})+\beta_{x_{1},f}y \quad{\rm
and}\quad\delta(x_{2}\otimes
f)y=\alpha_{x_{2},f}(y)(x_{2})+\beta_{x_{2},f}y.$$ By the additivity
of $\delta$, the above equations yield
$$\begin{array}{rl}0=&(\alpha_{x_{1}+x_{2},f}-\alpha_{x_{1},f})(y)x_{1}+(\alpha_{x_{1}+x_{2},f}-\alpha_{x_{1},f})(y)x_{2}\\
&+(\beta_{x_{1}+x_{2},f}-\beta_{x_{1},f}-\beta_{x_{2},f})y.\end{array}$$
Since $\ker(f)\cap N\supseteq N_-$ is infinite-dimensional, we can
choose $y\in \ker f\cap N$ such that $y\not\in{\rm
span}\{x_1,x_2\}$. It follows that
$\beta_{x_{1}+x_{2},f}=\beta_{x_{1},f}+\beta_{x_{2},f}$, that is,
$\beta_{x,f}$ is additive for $x\in{\mathcal D}_1({\mathcal N})$.
Now for any non-zero scalar $s$, by Eq.(2.14), we have
$$\delta(sx_1\otimes
f)y=s\alpha_{sx_1,f}(y)x_1+\beta_{sx_1,f}y\quad{\rm and}\quad
s\delta(x_1\otimes f)y=s\alpha_{x_1,f}(y)x_1+s\beta_{x_1,f}y,$$
which imply
$$s(\alpha_{sx_1,f}(y)-\alpha_{x_1,f}(y))x_1+(\beta_{sx_1,f}-s\beta_{x_1,f})y=0.$$
Still choosing $y$ linearly independent of $x_1$ gives
$\beta_{sx_1,f}=s\beta_{x_1,f}$. Hence $\beta_{x,f}$ is linear in
$x$.

For any $f_1,f_2\in{\mathcal D}_2({\mathcal N})$, there exists some
$N\in{\mathcal N}$ such that $f_{1},\ f_{2}\in N_{-}^{\perp}$. Take
any $x\in N$; and then $x\otimes f_{1}$, $x\otimes f_{2}\in {\rm
Alg}\mathcal {N}$. Since $f_1,f_2\in{\mathcal D}_2({\mathcal N})$,
$f_1,f_2\in M^\perp$ for some $M\in{\mathcal
N}\setminus\{\{0\},X\}$. Thus $M\subset\ker (f_1)\cap\ker(f_2)$ and
$\ker f_{1}\cap
 \ker f_{2}\cap N=M$ or $N$. So $\dim(\ker f_{1}\cap
 \ker f_{2}\cap N)=\infty$. For any $y\in \ker
f_{1}\cap
 \ker f_{2}\cap N$, by Eq.(2.14), one has
$$\delta(x\otimes
(f_{1}+f_{2}))y=\alpha_{x,f_{1}+f_{2}}(y)x+\beta_{x,f_{1}+f_{2}}y,$$
$$\delta(x\otimes
f_{1})y=\alpha_{x,f_{1}}(y)x+\beta_{x,f_{1}}y\quad{\rm and}\quad
\delta(x\otimes f_{2})y=\alpha_{x,f_{2}}(y)x+\beta_{x,f_{2}}y.$$
Comparing the above three equations, we get
$$0=(\alpha_{x,f_{1}+f_{2}}(y)-\alpha_{x,f_{1}}(y)-\alpha_{x,f_{2}}(y))x
+(\beta_{x,f_{1}+f_{2}}-\beta_{x,f_{1}}-\beta_{x,f_{2}})y.\eqno(2.15)$$
Choosing   $y$  linearly independent of $x$  in Eq.(2.15) entails
$\beta_{x,f_{1}+f_{2}}=\beta_{x,f_{1}}+\beta_{x,f_{2}}.$

For $x\otimes(sf_{1})$, by Eq.(2.14), we get $$\delta(x\otimes
sf_{1})y=\alpha_{x,sf_{1}}(y)x+\beta_{x,sf_{1}}y\quad{\rm and}\quad
s\delta(x_{1}\otimes f)y=s\alpha_{x,f_{1}}(y)x+s\beta_{x,f_{1}}y.$$
Thus
$$(\alpha_{x,sf_{1}}-s\alpha_{x,f_{1}})(y)x+(\beta_{x,sf_{1}}-s\beta_{x,f_{1}})y=0.$$
Since, for $x\in \ker f_{1}\cap
 \ker f_{2}\cap N$, we can find $y\in \ker f_{1}\cap
 \ker f_{2}\cap N$ so that $y$ is linearly independent of
 $x$, it follows that $\beta_{x,sf_{1}}=s\beta_{x,f_{1}}$. Hence $\beta_{x,f}$ is
 linear in $f$, completing the proof of (1).

Now let us prove the conclusion (2). It follows from (1) that, for
any $x\otimes f\in {\rm Alg}{\mathcal N}$, there exists a functional
$\alpha_{x,f}: \ker (f)\rightarrow{\mathbb F}$ such that
$\delta(x\otimes f)y-\beta_{x,f}y=\alpha_{x,f}(y)x$ for any $y\in
\ker(f)$. Similar to the proof of Lemma 2.2(2), we can find linear
maps $B: \mathcal D_1(\mathcal N)\rightarrow \mathcal D_1(\mathcal
N)$ and $C: \mathcal D_{2}(\mathcal N)\rightarrow \mathcal
D_{2}(\mathcal N)$ such that $\delta(x\otimes
f)-\beta_{x,f}I=x\otimes Cf+Bx\otimes f$ holds for all $x\otimes
f\in {\mbox {\rm Alg}}\mathcal N$.

We complete the proof of the lemma. \hfill$\Box$

The following  Lemmas 2.5-2.6 come from \cite{HQ1}. For the
completeness,  we give their proofs here.

{\bf Lemma 2.5.} {\it Let $\mathcal N$  be a nest on a  real or
complex Banach space $X$. Suppose that $\delta: \mbox{\rm
Alg}{\mathcal N}\rightarrow \mbox{\rm Alg}{\mathcal N}$ is a linear
map. If there exists an injective operator or an operator with dense
range $Z \in {\rm Alg}{\mathcal N}$ such that $\delta$ is derivable
at $Z$, then, $\delta(I)=0$.}

{\bf Proof.} Since $\delta$ is derivable at $Z$ and $Z=IZ=ZI$, we
have $\delta(Z)=\delta(I)Z+I\delta(Z)=\delta(Z)+Z\delta(I)$. So
$\delta(I)Z=Z\delta(I)=0$. If $Z$ is injective, by $Z\delta(I)=0$,
we get $\delta(I)=0$; if $Z$ is an operator with dense range, then,
by $\delta(I)Z=0$, we get again $\delta(I)=0$. \hfill$\Box$

{\bf Lemma 2.6.} {\it Let $\mathcal N$  be a nest on a complex
Banach space $X $ and $\delta: \mbox{\rm Alg}{\mathcal N}\rightarrow
\mbox{\rm Alg}{\mathcal N}$ be a linear map derivable at $Z\in{\rm
Alg}{\mathcal N}$.

If $Z$ is an operator with dense range, then}

(1) {\it for every idempotent operator $P\in \mbox{\rm Alg}{\mathcal
N},$ we have $\delta(PZ)=\delta(P)Z+P\delta(Z)$, and moreover,
$\delta(P)=\delta(P)P+P\delta(P)$;}

(2) {\it for every operator $N\in \mbox{\rm Alg}{\mathcal N}$ with
$N^2=0$, we have $\delta(NZ)=\delta(N)Z+N\delta(Z)$, and moreover,
$\delta(N)N+N\delta(N)=0$.}

{\it If $Z$ is an injective operator, then }

($1'$) {\it for every idempotent operator $P \in \mbox{\rm
Alg}{\mathcal N}$, we have $\delta(ZP)=\delta(Z)P+Z\delta(P)$, and
moreover, $\delta(P)=\delta(P)P+P\delta(P)$;}

($2'$) {\it for every operator $N\in \mbox{\rm Alg}{\mathcal N}$
with $N^2=0$, we have $\delta(ZN)=\delta(Z)N+Z\delta(N)$, and
moreover, $\delta(N)N+N\delta(N)=0$.}

{\bf Proof.} Let $P\in \mbox{Alg}{\mathcal N}$ be any idempotent
operator. If $Z$ is an operator with dense range, then by Lemma 2.5,
we have
$$\begin{array}{rl}
\delta(Z)=&\delta(I-\frac{1-\sqrt{3}i}{2}P)(I-\frac{1+\sqrt{3}i}{2}P)Z\\
&+(I-\frac{1-\sqrt{3}i}{2}P)\delta(Z-\frac{1+\sqrt{3}i}{2}PZ)\\
=&-\frac{1-\sqrt{3}i}{2}\delta(P)Z+\delta(P)PZ+\delta(Z)\\
&-\frac{1+\sqrt{3}i}{2}\delta(PZ)-\frac{1-\sqrt{3}i}{2}P\delta(Z)+P\delta(PZ)
\end{array}$$
since $Z=(I-\frac{1-\sqrt{3}i}{2}P)(I-\frac{1+\sqrt{3}i}{2}P)Z$.
Thus we get
$$-\frac{1-\sqrt{3}i}{2}\delta(P)Z+\delta(P)PZ
-\frac{1+\sqrt{3}i}{2}\delta(PZ)-\frac{1-\sqrt{3}i}{2}P\delta(Z)+P\delta(PZ)=0.\eqno(2.16)$$
On the other hand,
$Z=(I-\frac{1+\sqrt{3}i}{2}P)(I-\frac{1-\sqrt{3}i}{2}P)Z$ gives
$$-\frac{1+\sqrt{3}i}{2}\delta(P)Z+\delta(P)PZ
-\frac{1-\sqrt{3}i}{2}\delta(PZ)-\frac{1+\sqrt{3}i}{2}P\delta(Z)+P\delta(PZ)=0.\eqno(2.17)$$
Combining Eq.(2.16) with Eq.(2.17), we get
$$\delta(PZ)=\delta(P)Z+P\delta(Z).$$
Replacing $\delta(PZ)$ by $\delta(P)Z+P\delta(Z)$ in Eq.(2.16), one
obtains $\delta(P)Z=\delta(P)PZ+P\delta(P)Z$. Note that $Z$ is an
operator with dense range. It follows that
$\delta(P)=\delta(P)P+P\delta(P).$ This completes the proof of
assertion (1).

If $Z$ is an injective operator, then by the equation
$$Z=Z(I-\frac{1-\sqrt{3}i}{2}P)(I-\frac{1+\sqrt{3}i}{2}P)=Z(I-\frac{1+\sqrt{3}i}{2}P)(I-\frac{1-\sqrt{3}i}{2}P),$$
using a similar argument as  above, one can get that
$\delta(ZP)=\delta(Z)P+Z\delta(P)$ and
$\delta(P)=\delta(P)P+P\delta(P).$ Hence  (1$'$) holds true.

For every operator $N\in \mbox{\rm Alg}{\mathcal N}$ with $N^2=0$,
if $Z$ is an operator with dense range, then, noting that
$Z=(I-N)(I+N)Z=(I+N)(I-N)Z$, we have
$$\delta(N)Z-\delta(N)NZ-\delta(NZ)+N\delta(Z)-N\delta(NZ)=0\eqno(2.18)$$
and
$$-\delta(N)Z-\delta(N)NZ+\delta(NZ)-N\delta(Z)-N\delta(NZ)=0$$
since $\delta$ is derivable at $Z$. Comparing the above two
equations, one gets$$\delta(NZ)=\delta(N)Z+N\delta(Z).$$ Replacing
$\delta(NZ)$ by $\delta(N)Z+N\delta(Z)$ in Eq.(2.18) and noting that
the range of $Z$ is dense, it follows that
$\delta(N)N+N\delta(N)=0$.

If $Z$ is an injective operator, then by the equation
$Z=Z(I-N)(I+N)=Z(I+N)(I-N)$, a similar argument shows that
$\delta(ZN)=\delta(Z)N+Z\delta(N)$ and
$\delta(P)=\delta(P)P+P\delta(P).$ Hence the assertions  (2) and
(2$'$) hold true,
 completing the proof. \hfill$\Box$

By Lemmas 2.5-2.6,  Lemmas 2.2-2.4 hold for linear maps between nest
algebras on complex Banach space derivable at an injective operator
or an operator with dense range.

\section{Proof of Theorem 1.1}

In this section we complete the proof of Theorem 1.1.

{\bf Proof of Theorem 1.1.} The ``if " part is obvious. We only need
to prove the ``only if" part.

In the following, we always assume that $Z\in{\rm Alg}{\mathcal N}$
is an injective operator or an operator with dense range and
$\delta:{\rm Alg}{\mathcal N}\rightarrow{\rm Alg}{\mathcal N}$ is a
linear map derivable at $Z$. Then, for any invertible element $A\in
{\rm Alg}\mathcal N$, since $Z=AA^{-1}Z=ZA^{-1}A$ and $\delta$ is
derivable in $Z$, we have
$$\delta(Z)=\delta(A)A^{-1}Z+A\delta(A^{-1}Z)$$
and $$\delta(Z)=\delta(ZA^{-1})A+ZA^{-1}\delta(Z),$$ which imply
respectively that
 $$\delta (A^{-1}Z)=A^{-1}\delta(Z)-A^{-1}\delta
(A)A^{-1}Z \eqno(3.1)$$ and
$$\delta (ZA^{-1})=\delta(Z)A^{-1}-ZA^{-1}\delta
(A)A^{-1}.\eqno(3.2)$$

In the sequel, we will prove that $\delta$ is a derivation by
considering three cases.

{\bf Case 1}. $X_-\not= X.$

In this case, by Lemmas 2.2 and 2.6, there exist  linear
transformations $B:X\rightarrow X$ and  $C:X_{-}^{\perp}\rightarrow
X^{\ast}$ such that $\delta(x\otimes f)=Bx\otimes f+x\otimes Cf $
and $\langle Bx,f \rangle +\langle x,Cf\rangle=0 $ for all $x\in X$
and $f\in X_{-}^{\perp}$. Thus, by Lemma 2.6 and the linearity of
$\delta$, if $Z$ is an operator with dense range,  we have
$$\delta(x\otimes fZ)=(Bx\otimes f)Z+(x\otimes Cf)Z+(x\otimes
f)\delta(Z); \eqno(3.3)$$ if $Z$ is an injective operator, we have
$$\delta(Zx\otimes f)=\delta(Z)(x\otimes f)+ZBx\otimes f+Zx\otimes
Cf. \eqno(3.4)$$

Note that, for any $T$ and any $x\otimes f\in{\rm Alg} {\mathcal
N}$, there exists some $\lambda \in{\mathbb C}$ such that
$|\lambda|>\|T\|$ and $\|(\lambda I-T)^{-1}x\|\|f\|<1$. Then both
$\lambda I-T$ and $\lambda I-T-x\otimes f=(\lambda I-T)(I-(\lambda
I-T)^{-1}x\otimes f)$ are invertible with their inverses are still
in ${\rm Alg}{\mathcal N}$. It is obvious that $(I-(\lambda
I-T)^{-1}x\otimes f)^{-1}=I+(1-\alpha )^{-1}(\lambda
I-T)^{-1}x\otimes f$, where $\alpha=\langle (\lambda
I-T)^{-1}x,f\rangle$.

{\bf Claim 1.1.} $\delta$ is a derivation if $Z$ is an operator with
dense range.

For any $T$ and any $x\otimes f\in{\rm Alg} {\mathcal N}$ with $x\in
X$ and $f\in X_-^\perp$, take $\lambda \in{\mathbb C}$ such that
$|\lambda|>\|T\|$ and $\|(\lambda I-T)^{-1}x\|\|f\|<1$. By Lemma
2.2, Eqs.(3.1), (3.3)  and the fact $\delta(I)=0$ (Lemma 2.5), we
have
$$\begin{array}{rl} \delta(Z)=&\delta(\lambda I-T-x\otimes f)(I+(1-\alpha
)^{-1}(\lambda I-T)^{-1}x\otimes f)(\lambda I-T)^{-1}Z\\
& +(\lambda I-T-x\otimes f)\delta ((I+(1-\alpha )^{-1}(\lambda
I-T)^{-1}x\otimes f)(\lambda I-T)^{-1}Z)\\
=&[-\delta (T)-Bx\otimes f-x\otimes Cf][(\lambda
I-T)^{-1}Z\\
&+(1-\alpha)^{-1}(\lambda I-T)^{-1}(x\otimes f)(\lambda
I-T)^{-1}Z]\\
&+(\lambda I-T-x\otimes f)[(\lambda I-T)^{-1}\delta(Z)+(\lambda
I-T)^{-1}\delta(T)(\lambda I-T)^{-1}Z\\
&+(1-\alpha)^{-1}B(\lambda I-T)^{-1}(x\otimes f)(\lambda
I-T)^{-1}Z\\
&+(1-\alpha)^{-1}(\lambda I-T)^{-1}x\otimes C((\lambda
I-T^*)^{-1}f)Z\\
&+(1-\alpha)^{-1}(\lambda I-T)^{-1}x\otimes f(\lambda
I-T)^{-1}\delta(Z)]\\ =&\delta(Z)-(1-\alpha)^{-1}B(x\otimes
(\lambda I-T^*)^{-1}f)Z\\
&-(1-\alpha)^{-1}\delta(T)(\lambda I-T)^{-1}(x\otimes (\lambda
I-T^*)^{-1}f)Z\\
&-(x\otimes (\lambda I-T^*)^{-1}Cf)Z+(x\otimes C(\lambda
I-T^*)^{-1}f)Z\\
&+(1-\alpha)^{-1}(\lambda I-T)B(\lambda I-T)^{-1}(x\otimes
(\lambda I-T^*)^{-1}f)Z\\
&-(x\otimes (\lambda I-T^*)^{-1}\delta(T)^*(\lambda I-T^*)^{-1}f)Z\\
&-(1-\alpha)^{-1}(\langle(\lambda I-T)^{-1}x, Cf\rangle+\langle
B(\lambda I-T)^{-1}x,f\rangle)(x\otimes (\lambda I-T^*)^{-1}f)Z .
\end{array}$$
As $\langle(\lambda I-T)^{-1}x, Cf\rangle+\langle B(\lambda
I-T)^{-1}x,f\rangle=0$, the above equation becomes
$$\begin{array}{rl} 0=&(1-\alpha)^{-1}B(x\otimes
(\lambda I-T^*)^{-1}f)Z+(1-\alpha)^{-1}\delta(T)(\lambda
I-T)^{-1}(x\otimes
(\lambda I-T^*)^{-1}f)Z\\
&+(x\otimes (\lambda I-T^*)^{-1}Cf)Z-(x\otimes
C(\lambda I-T^*)^{-1}f)Z\\
&-(1-\alpha)^{-1}(\lambda I-T)B(\lambda I-T)^{-1}(x\otimes
(\lambda I-T^*)^{-1}f)Z\\
&+(x\otimes (\lambda I-T^*)^{-1}\delta(T)^*(\lambda I-T^*)^{-1}f)Z.
\end{array}
$$
Since the range of $Z$ is dense, it follows  that
$$\begin{array}{rl} 0=&(1-\alpha)^{-1}Bx\otimes
(\lambda I-T^*)^{-1}f+(1-\alpha)^{-1}\delta(T)(\lambda
I-T)^{-1}x\otimes
(\lambda I-T^*)^{-1}f\\
&+x\otimes (\lambda I-T^*)^{-1}Cf-x\otimes
C(\lambda I-T^*)^{-1}f\\
&-(1-\alpha)^{-1}(\lambda I-T)B(\lambda I-T)^{-1}x\otimes
(\lambda I-T^*)^{-1}f\\
&+x\otimes (\lambda I-T^*)^{-1}\delta(T)^*(\lambda I-T^*)^{-1}f,
\end{array}
$$
and so
$$\begin{array}{rl}
&[\delta(T)(\lambda I-T)^{-1}-(\lambda I-T)B(\lambda
I-T)^{-1}+B]x\otimes (\lambda
I-T^*)^{-1}f \\
=&x\otimes(1-\alpha)[C(\lambda I-T^*)^{-1}-(\lambda
I-T^*)^{-1}C-(\lambda I-T^*)^{-1}\delta(T)^*(\lambda I-T^*)^{-1}]f.
\end{array}
$$
Hence $[\delta(T)(\lambda I-T)^{-1}-(\lambda I-T)B(\lambda
I-T)^{-1}+B]x$ is linearly dependent of $x$ for every $x\in X$. This
entails that there is a scalar $\beta$ such that
$$\delta(T)(\lambda I-T)^{-1}-(\lambda I-T)B(\lambda
I-T)^{-1}+B=\beta I $$ on $X$. It follows that
$\delta(T)=BT-TB+\beta (\lambda I-T).$ By taking different $\lambda$
in the equation, we see that $\beta=0$ and consequently $\delta(T)
=BT -TB$ holds for all $T\in {\rm Alg}{\mathcal N}$, that is,
$\delta$ is a derivation.

{\bf Claim 1.2.} $\delta$ is a derivation if $Z$ is an injective
operator.

For any $T$ and any $x\otimes f\in{\rm Alg} {\mathcal N}$ with $x\in
X$ and $f\in X_-^\perp$, take $\lambda \in{\mathbb C}$ such that
$|\lambda|>\|T\|$ and $\|(\lambda I-T)^{-1}x\|\|f\|<1$. Note that
$(\lambda I-T^*)^{-1}f\in X_-^\perp$ and
 $Z=Z(I+(1-\alpha
)^{-1}(\lambda I-T)^{-1}x\otimes f)(\lambda I-T)^{-1}(\lambda
I-T-x\otimes f)$. Then, by Lemma 2.2, Eqs.(3.2), (3.4) and the fact
$\delta(I)=0$, we have
$$\begin{array}{rl} \delta(Z)=&\delta(Z(\lambda I-T)^{-1}+(1-\alpha
)^{-1}Z(\lambda I-T)^{-1}x\otimes f(\lambda I-T)^{-1})(\lambda
I-T-x\otimes f)\\
& +(Z(\lambda I-T)^{-1}+(1-\alpha )^{-1}Z(\lambda I-T)^{-1}x\otimes
f(\lambda I-T)^{-1})\delta(\lambda I-T-x\otimes f)\\
=&[\delta(Z)(\lambda I-T)^{-1}+Z(\lambda I-T)^{-1}\delta(T)(\lambda
I-T)^{-1}\\
&+(1-\alpha )^{-1}\delta(Z)(\lambda I-T)^{-1}(x\otimes
f)(\lambda I-T)^{-1}\\
&+(1-\alpha )^{-1}ZB(\lambda I-T)^{-1}(x\otimes f)(\lambda
I-T)^{-1}\\
&+(1-\alpha )^{-1}Z(\lambda I-T)^{-1}(x\otimes C(\lambda
I-T^*)^{-1}f)][(\lambda
I-T)-x\otimes f]\\
&+[Z(\lambda I-T)^{-1}\\&+(1-\alpha )^{-1}Z(\lambda
I-T)^{-1}x\otimes
f(\lambda I-T)^{-1}][-\delta(T)-Bx\otimes f-x\otimes Cf]\\
=&\delta(Z)-Z(\lambda I-T)^{-1}\delta(T)(\lambda I-T)^{-1}x\otimes
f\\
&+(1-\alpha )^{-1}Z(\lambda I-T)^{-1}(x\otimes (\lambda
I-T^*)C(\lambda
I-T^*)^{-1}f)\\
&+ZB(\lambda I-T)^{-1}x\otimes f-Z(\lambda I-T)^{-1}Bx\otimes f\\
&-(1-\alpha )^{-1}Z(\lambda I-T)^{-1}x\otimes \delta(T)^*(\lambda
I-T^*)^{-1}f-(1-\alpha )^{-1}Z(\lambda I-T)^{-1}x\otimes Cf.
\end{array}$$
\if Here, in the last equality we apply Eq.(8) for rank one operator
$x\otimes (\lambda I-T^*)^{-1}f\in{\rm Alg}{\mathcal N}$.\fi It
follows  that
$$\begin{array}{rl} 0=&Z(\lambda I-T)^{-1}\delta(T)(\lambda I-T)^{-1}x\otimes
f\\
&-(1-\alpha )^{-1}Z(\lambda I-T)^{-1}x\otimes (\lambda
I-T^*)C(\lambda
I-T^*)^{-1}f\\
&-ZB(\lambda I-T)^{-1}x\otimes f+Z(\lambda I-T)^{-1}Bx\otimes f\\
&+(1-\alpha )^{-1}Z(\lambda I-T)^{-1}x\otimes \delta(T)^*(\lambda
I-T^*)^{-1}f+(1-\alpha )^{-1}Z(\lambda I-T)^{-1}x\otimes Cf.
\end{array}
$$
Since $\ker Z=\{0\}$, we get
$$\begin{array}{rl} 0=&(\lambda I-T)^{-1}\delta(T)(\lambda I-T)^{-1}x\otimes
f\\
&-(1-\alpha )^{-1}(\lambda I-T)^{-1}x\otimes (\lambda
I-T^*)C(\lambda
I-T^*)^{-1}f\\
&-B(\lambda I-T)^{-1}x\otimes f+(\lambda I-T)^{-1}Bx\otimes f\\
&+(1-\alpha )^{-1}(\lambda I-T)^{-1}x\otimes \delta(T)^*(\lambda
I-T^*)^{-1}f+(1-\alpha )^{-1}(\lambda I-T)^{-1}x\otimes Cf.
\end{array}
$$
Multiplying the above equation by $(\lambda I-T)$ from the left, one
has
$$\begin{array}{rl} 0=&\delta(T)(\lambda I-T)^{-1}x\otimes
f-(1-\alpha )^{-1}x\otimes (\lambda I-T^*)C(\lambda
I-T^*)^{-1}f\\
&-(\lambda I-T)B(\lambda I-T)^{-1}x\otimes f+Bx\otimes f\\
&+(1-\alpha )^{-1}x\otimes \delta(T)^*(\lambda
I-T^*)^{-1}f+(1-\alpha )^{-1}x\otimes Cf.
\end{array}$$ That is,
$$\begin{array}{rl}
&[\delta(T)(\lambda I-T)^{-1}-(\lambda I-T)B(\lambda I-T)^{-1}+B]x\otimes f \\
=&(1-\alpha )^{-1}x\otimes [(\lambda I-T^*)C(\lambda
I-T^*)^{-1}-\delta(T)^*(\lambda I-T^*)^{-1}-C]f.
\end{array}
$$
Now using the same argument as in the proof of Claim 1.1, we see
that $\delta$ is a derivation.

{\bf Case 2}. $\{0\}\not=\{0\}_+$.

 In this case, one can use Lemmas 2.3,  2.5-2.6 and a similar argument of
Case 1 to check that $\delta$ is a derivation, and we omit the
details here.

\if false

For any $T\in\mbox{\rm Alg}{\mathcal N}$ and $x\otimes f
\in\mbox{\rm Alg}{\mathcal N}$ with $f\in X^{\ast},x\in \{0\}_+$,
take  a scalar  $\lambda$ such that $|\lambda|>\|T\|$ and
$\|(\lambda I-T)^{-1}x\|\|f\|<1$. Still write $\alpha=\langle
(\lambda I-T)^{-1}x,f\rangle$. As $\delta$ is derivable at $I$, we
have
$$\begin{array}{rl} 0=&\delta(\lambda I-T-x\otimes f)(I+(1-\alpha
)^{-1}(\lambda I-T)^{-1}x\otimes f)(\lambda I-T)^{-1}\\
& +(\lambda I-T-x\otimes f)\delta ((I+(1-\alpha )^{-1}(\lambda
I-T)^{-1}x\otimes f)(\lambda I-T)^{-1}).
\end{array}$$ Taking conjugate gives
$$\begin{array}{rl} 0=&-(1-\alpha)^{-1}(\lambda
I-T^{\ast})^{-1}f\otimes ((\lambda
I-T)^{-1}x)^{\ast\ast}\delta(T)^{\ast} -(\lambda
I-T^{\ast})^{-1}\delta(x\otimes f)^{\ast}\\
&-(1-\alpha)^{-1}(\lambda I-T^{\ast})^{-1}f\otimes ((\lambda
I-T)^{-1}x)^{\ast\ast}\delta(x\otimes f)^{\ast}\\
&+(1-\alpha)^{-1}\delta((\lambda I-T)^{-1})x\otimes(\lambda
I-T^{\ast})^{-1}f )^{\ast}(\lambda I-T^{\ast})\\
&-((\lambda I-T)^{-1})\delta(T)(\lambda I-T)^{-1})^{\ast}f \otimes
x^{\ast\ast}\\
&-(1-\alpha)^{-1}\delta((\lambda I-T)^{-1})^{-1}x\otimes f (\lambda
I-T)^{-1}))^{\ast}(f\otimes x^{\ast\ast}).
\end{array}$$
Applying Claim 2.2 to the above equation, one  gets
$$\begin{array}{rl} 0=&-(1-\alpha)^{-1}(\lambda
I-T^{\ast})^{-1}f\otimes ((\lambda
I-T)^{-1}x)^{\ast\ast}\delta(T)^{\ast} \\
&-(\lambda I-T^{\ast})^{-1}Bf\otimes x^{\ast\ast}-(\lambda
I-T^{\ast})^{-1}f\otimes Cx^{\ast\ast}\\
&-(1-\alpha)^{-1}\langle Bf,((\lambda
I-T)^{-1}x)^{\ast\ast}\rangle(\lambda I-T^{\ast})^{-1}f\otimes
x^{\ast\ast}\\
&-(1-\alpha)^{-1}\langle f,((\lambda I-T)^{-1}x)^{\ast\ast}\rangle
(\lambda
I-T^{\ast})^{-1}f\otimes Cx^{\ast\ast}\\
&+(1-\alpha)^{-1}B((\lambda I-T)^{-1})f\otimes((\lambda
I-T)^{-1}x)^{\ast\ast}(\lambda I-T^{\ast})\\
&+(1-\alpha)^{-1}(\lambda I-T^{\ast})^{-1})f\otimes C(\lambda
I-T)^{-1}x)^{\ast\ast}(\lambda I-T^{\ast})\\
&-(\lambda I-T^{\ast})^{-1}\delta(T)^{\ast}(\lambda
I-T^{\ast})^{-1}f\otimes x^{\ast\ast}\\
&-(1-\alpha)^{-1}\langle f,((\lambda
I-T)^{-1}x)^{\ast\ast}\rangle B((\lambda I-T)^{-1})f\otimes x^{\ast\ast}\\
&-(1-\alpha)^{-1}\langle f,C(\lambda I-T)^{-1}x)^{\ast\ast}\rangle
(\lambda I-T^{\ast})^{-1}f\otimes x^{\ast\ast}.
\end{array}$$
Applying Claim 2.3, the equation above becomes
$$\begin{array}{rl} 0=&-(1-\alpha)^{-1}(\lambda
I-T^{\ast})^{-1}f\otimes ((\lambda
I-T)^{-1}x)^{\ast\ast}\delta(T)^{\ast}\\
&-(\lambda I-T^{\ast})^{-1}Bf\otimes
x^{\ast\ast}-(1-\alpha)^{-1}(\lambda
I-T^{\ast})^{-1}f\otimes Cx^{\ast\ast}\\
&+ B((\lambda I-T)^{-1}f)\otimes x^{\ast\ast}-(\lambda
I-T^{\ast})^{-1}\delta(T)^{\ast}(\lambda
I-T^{\ast})^{-1}f\otimes x^{\ast\ast}\\
&+(1-\alpha)^{-1}(\lambda I-T^{\ast})^{-1})f\otimes C(\lambda
I-T)^{-1}x)^{\ast\ast}(\lambda I-T^{\ast}).
\end{array}$$
Hence,
$$\begin{array}{rl}0=&(\lambda
I-T^{\ast})^{-1}f\otimes [(\lambda
I-T)^{-1}x)^{\ast\ast}\delta(T)^{\ast}+Tx^{\ast\ast}-C(\lambda
I-T)^{-1}x)^{\ast\ast}(\lambda I-T^{\ast})]\\
&+(1-\alpha)[(\lambda I-T^{\ast})^{-1}Bf- B((\lambda
I-T)^{-1}f)+(\lambda I-T^{\ast})^{-1}\delta(T)^{\ast}(\lambda
I-T^{\ast})^{-1}f]\otimes x^{\ast\ast}.
\end{array}$$
This implies that $[(\lambda I-T^{\ast})^{-1}B- B(\lambda
I-T)^{-1}+(\lambda I-T^{\ast})^{-1}\delta(T)^{\ast}(\lambda
I-T^{\ast})^{-1}]f $ is linearly dependent of $(\lambda
I-T^{\ast})^{-1}f$ for every $f\in X^*$. So there exists a scalar
$p_{\lambda}$ such that
$$[(\lambda I-T^{\ast})^{-1}B- B(\lambda
I-T)^{-1}+(\lambda I-T^{\ast})^{-1}\delta(T)^{\ast}(\lambda
I-T^{\ast})^{-1}]f=p_{\lambda}(\lambda I-T^{\ast})^{-1}f$$ for all
$f\in X^{\ast}$. It follows that
$\delta(T)^{\ast}=T^{\ast}B-BT^{\ast}+(1-\alpha)^{-1}p_{\lambda}(\lambda
I-T^{\ast})^{-1}I.$ By taking different $\lambda$, we see that
$p_{\lambda}=0$ and consequently,
$$\delta(T)^{\ast}=T^{\ast}B-BT^{\ast}\quad{\rm holds\ \ for\ \ all}\quad T\in\mbox{\rm Alg}{\mathcal
N}. \eqno(3.5)$$

Now for  any $T,S\in \mbox{\rm Alg}{\mathcal N}$, by Eq.(3.5),  we
have
$$\delta(TS)^{\ast}=S^{\ast}T^{\ast}B-BT^{\ast}S^{\ast}$$ and
$$(\delta(T)S+T\delta(S))^{\ast}=S^{\ast}(T^{\ast}B-BT^{\ast})+(S^{\ast}B-BS^{\ast})T^{\ast}
=S^{\ast}T^{\ast}B-BT^{\ast}S^{\ast}.$$ Comparing the above two
equations, one gets
$\delta(TS)^{\ast}=(\delta(T)S+T\delta(S))^{\ast}$, and so
$\delta(TS)=\delta(T)S+T\delta(S)$   holds for all $T,S\in \mbox{\rm
Alg}{\mathcal N}$.

Therefore for the case $\{0\}\not=\{0 \}_{+}$, $\delta$ is still a
derivation.\fi

We remark that, the trivial case, that is, Alg${\mathcal
N}={\mathcal B}(X)$, is included in both Case 1 and Case 2.

Finally, let us consider the case that both end points of the nest
are limit points.

{\bf Case 3}. $\{0\}=\{0\}_+$ and $X_-=X$.

In this case  $X$ is infinite dimensional and every nonzero element
$N$ in $\mathcal N$ is infinite dimensional.

By Lemmas 2.4 and 2.6, there exists a bilinear functional $\beta:
({\mathcal D}_1({\mathcal N})\times {\mathcal D}_2({\mathcal
N}))\cap{\rm Alg}{\mathcal N}\rightarrow \mathbb C$, linear
transformations $B: \mathcal D_1(\mathcal N)\rightarrow \mathcal
D_1(\mathcal N)$ and $C: \mathcal D_{2}(\mathcal N)\rightarrow
\mathcal D_{2}(\mathcal N)$ such that $$(\delta(x\otimes
f)-\beta_{x,f}I)\ker (f)\subseteq \mbox{\rm span}\{x\}$$ and
$$\delta(x\otimes f)-\beta_{x,f}I=x\otimes Cf+Bx\otimes f\eqno(3.5)$$
hold for all $x\otimes f\in {\mbox {\rm Alg}}\mathcal N$.  Then, by
Lemma 2.6 and the linearity of $\delta$, if $Z$ is an operator with
dense range,  we have
$$\delta(x\otimes fZ)=(Bx\otimes f)Z+(x\otimes Cf)+\beta_{x,f}Z+(x\otimes
f)\delta(Z); \eqno(3.6)$$ if $Z$ is an injective operator,  we have
$$\delta(Zx\otimes f)=\delta(Z)(x\otimes f)+ZBx\otimes f+Zx\otimes
Cf+\beta_{x,f}Z. \eqno(3.7)$$

Now let $\tilde{\beta}_{x,f}=\langle x,Cf\rangle+\langle
Bx,f\rangle$. We claim that
$$\left\{
\begin{array}{ll}
  \tilde{\beta}_{x,f}=0, & \mbox{\rm if}\ \langle x, f\rangle\not=0;  \\
  \tilde{\beta}_{x,f}=-2\beta_{x,f}, & \mbox{\rm if}\ \langle x, f\rangle=0.\end{array} \right. \eqno(3.8)$$
  In fact, if $\langle x, f\rangle\not=0$,
then, by Lemma 2.6(1) and (1$'$),
 $(x\otimes f)\delta (x\otimes f)(x\otimes
f)=0$, and hence $\langle x,Cf\rangle+\langle Bx,f
\rangle=-\beta_{x,f}=0$; if $\langle x, f\rangle=0$, then, by Lemma
2.6(2) and (2$'$),   $\delta(x\otimes f)(x\otimes f)+(x\otimes
f)\delta(x\otimes f)=0$, which, together with Eq.(3.5), implies that
$\langle x,Cf\rangle+\langle Bx,f\rangle=-2\beta_{x,f}$.

For any $T$ and $x\otimes f\in {\rm Alg}{\mathcal N}$, take
$\lambda$ such that $|\lambda|>\|T\|$ and $\|(\lambda
I-T)^{-1}x\|\|f\|<1$. \if Then both $\lambda I-T$ and $\lambda
I-T-x\otimes f=(\lambda I-T)(I-(\lambda I-T)^{-1}x\otimes f)$ are
invertible with their inverses   still in ${\rm Alg}{\mathcal
N}$.\fi Note that
$$I=(\lambda I-T-x\otimes f)(I+(1-\alpha )^{-1}(\lambda
I-T)^{-1}x\otimes f)(\lambda I-T)^{-1},$$ where $\alpha=\langle
(\lambda I-T)^{-1}x,f\rangle$.

We prove that $\delta$ is a derivation respectively by assume that
$Z$ is injective or of dense range.

{\bf Claim 3.1.}   $\delta$ is a derivation if $Z$ has dense range.

 Since $\delta$ is derivable at   $Z$ of dense range, by Lemma 2.5, Eqs.(3.1) and (3.3), we have
$$\begin{array}{rl} \delta(Z)=&\delta(\lambda I-T-x\otimes f)[I+(1-\alpha
)^{-1}(\lambda I-T)^{-1}x\otimes f](\lambda I-T)^{-1}Z\\
& +(\lambda I-T-x\otimes f)\delta ((I+(1-\alpha )^{-1}(\lambda
I-T)^{-1}x\otimes f)(\lambda I-T)^{-1}Z)\\
=&[-\delta (T)-Bx\otimes f-x\otimes Cf-\beta_{x,f}I][(\lambda
I-T)^{-1}Z\\&+(1-\alpha)^{-1}(\lambda I-T)^{-1}(x\otimes f)(\lambda
I-T)^{-1}Z]\\
&+[\lambda I-T-x\otimes f][(\lambda I-T)^{-1}\delta(Z)+(\lambda
I-T)^{-1}\delta(T)(\lambda
I-T)^{-1} Z\\
&+(1-\alpha)^{-1}(B(\lambda I-T)^{-1}(x\otimes f)(\lambda
I-T)^{-1}Z\\
&+(\lambda I-T)^{-1}x\otimes ((\lambda I-T)^{-1})^*f
\delta(Z)+\beta_{(\lambda I-T)^{-1}x,((\lambda I-T)^{-1})^*f}Z)].
\end{array}$$
Note that $(1-\alpha )^{-1}\alpha =(1-\alpha )^{-1}-1$ and
$\tilde{\beta}_{x,f}=\langle x,Cf\rangle+\langle Bx,f\rangle$. It
follows  that
$$\begin{array}{rl}&\beta_{x,f}(\lambda I-T)^{-1}Z-(1-\alpha)^{-1}\beta_{(\lambda I-T)^{-1}x,(\lambda
I-T^*)^{-1}f}(\lambda I-T)Z\\
=&(1-\alpha)^{-1}\delta(\lambda I-T)((\lambda I-T)^{-1}x\otimes
(\lambda I-T^*)^{-1}f)Z\\
&-(1-\alpha)^{-1}(Bx\otimes (\lambda
I-T^*)^{-1}f)Z-(x\otimes(\lambda
I-T^*)^{-1}Cf)Z\\
&-(1-\alpha)^{-1}\tilde{\beta}_{(\lambda I-T)^{-1}x,f}
(x\otimes(\lambda I-T^*)^{-1}f)Z\\
&-(1-\alpha)^{-1}\beta_{x,f}((\lambda I-T)^{-1}x\otimes (\lambda
I-T^*)^{-1}f)Z\\
&+(x\otimes f(\lambda I-T)^{-1}\delta(\lambda I-T)(\lambda
I-T)^{-1})Z+(x\otimes C(\lambda I-T^*)^{-1}f)Z\\
&+(1-\alpha)^{-1}((\lambda I-T)B(\lambda I-T)^{-1}x\otimes (\lambda
I-T^*)^{-1}f)Z\\
&-(1-\alpha)^{-1}\beta_{(\lambda I-T)^{-1}x,(\lambda
I-T^*)^{-1}f}(x\otimes f)Z.
\end{array} $$
Since the range of $Z$ is dense, the above equation implies
$$\begin{array}{rl}&\beta_{x,f}(\lambda I-T)^{-1}-(1-\alpha)^{-1}\beta_{(\lambda I-T)^{-1}x,(\lambda
I-T^*)^{-1}f}(\lambda I-T)\\
=&(1-\alpha)^{-1}\delta(\lambda I-T)(\lambda I-T)^{-1}x\otimes
(\lambda I-T^*)^{-1}f\\
&-(1-\alpha)^{-1}Bx\otimes (\lambda I-T^*)^{-1}f-x\otimes(\lambda
I-T^*)^{-1}Cf\\
&-(1-\alpha)^{-1}\tilde{\beta}_{(\lambda I-T)^{-1}x,f}
x\otimes(\lambda I-T^*)^{-1}f\\
&-(1-\alpha)^{-1}\beta_{x,f}(\lambda I-T)^{-1}x\otimes (\lambda
I-T^*)^{-1}f\\
&+x\otimes f(\lambda I-T)^{-1}\delta(\lambda I-T)(\lambda
I-T)^{-1}+x\otimes C(\lambda I-T^*)^{-1}f\\
&+(1-\alpha)^{-1}(\lambda I-T)B(\lambda I-T)^{-1}x\otimes (\lambda
I-T^*)^{-1}f\\
&-(1-\alpha)^{-1}\beta_{(\lambda I-T)^{-1}x,(\lambda
I-T^*)^{-1}f}(x\otimes f).
\end{array}\eqno (3.9)$$

Next we show that $\beta_{x,f}=0$ for any  $x\otimes f\in{\rm Alg}
{\mathcal N}$.    If $\langle x,f\rangle\neq0$, it is obvious that
$\beta_{x,f}=0$. For the case $\langle x,f\rangle=0$, to prove
$\beta_{x,f}=0$, we simplify Eq.(3.9) by letting $W=(\lambda
I-T)^{-1}$. Then $\alpha=\langle Wx,f\rangle$ and Eq.(3.9) becomes
$$\begin{array}{rl}&(1-\langle Wx,f\rangle)(\beta_{x,f}W+x\otimes
W^{*}Cf-x\otimes fW\delta(W^{-1})W-x\otimes CW^{*}f)\\
=&\beta_{Wx,W^{*}f}W^{-1}+\delta(W^{-1})Wx\otimes W^{*}f-Bx\otimes
W^{*}f-\tilde{\beta}_{Wx,f}x\otimes W^{*}f\\
&-\beta_{x,f}Wx\otimes W^{*}f+W^{-1}BW^{-1}x\otimes
W^{*}f-\beta_{Wx,W^{*}f}x\otimes f.
\end{array}\eqno (3.10)$$
Let $t$ be a real or complex number with $t\neq 0,1$.  Replacing $x$
by $tx$ in Eq.(3.10) and using the bilinearity  of $\beta$, we get
$$\begin{array}{rl}&t(1-t\langle
Wx,f\rangle)(\beta_{x,f}W+x\otimes W^{*}Cf-x\otimes
fW\delta(W^{-1})W-x\otimes CW^{*}f)\\
=&t\beta_{Wx,W^{*}f}W^{-1}+t\delta(W^{-1})Wx\otimes
W^{*}f-tBx\otimes
W^{*}f-t^2\tilde{\beta}_{Wx,f}x\otimes W^{*}f\\
&-t^2\beta_{x,f}Wx\otimes W^{*}f+tW^{-1}BW^{-1}x\otimes
W^{*}f-t^2\beta_{Wx,W^{*}f}x\otimes f.
\end{array} \eqno (3.11)$$
Comparing Eq.(3.10) and Eq.(3.11) gives
$$\begin{array}{rl}\langle Wx,f\rangle \beta_{x,f}W=&
-\langle Wx,f\rangle x\otimes W^{*}Cf +\langle Wx,f\rangle x\otimes
fW\delta(W^{-1})W\\
&+\langle Wx,f\rangle x\otimes CW^{*}f+\tilde{\beta}_{Wx,f}x\otimes
W^{*}f\\
&+\beta_{x,f}Wx\otimes W^{*}f-\beta_{Wx,W^{*}f}x\otimes f.
\end{array} \eqno (3.12)$$
Note that, the right side of Eq.(3.12) is a finite rank operator.
Thus, we must have $\langle Wx,f\rangle \beta_{x,f}=0$. If $\langle
Wx,f\rangle\neq0$, then $\beta_{x,f}=0$; if $\langle Wx,f\rangle=0$,
then   Eq.(3.12) leads to
$$0=-2{\beta}_{Wx,f}x\otimes W^*f+\beta_{x,f}Wx\otimes
W^*f-\beta_{Wx,W^{*}f}x\otimes  f.\eqno (3.13)$$  Assume on the
contrary that $\beta_{x,f}\neq0$. Then Eq.(3.13) gives
$$Wx\otimes W^*f=x\otimes(\frac{2\beta_{Wx,f}}{\beta_{x,f}}W^*f+\frac{\beta_{Wx,W^{*}f}}{\beta_{x,f}}f), $$
which implies that, for any $x\otimes f\in$Alg$\mathcal N$ with
$x\in  \ker f$, there exists a scalar $t_{x,f}\not=0$ such that
$Wx=t_{x,f}x$. It follows that $Tx= \xi_{x,f,\lambda}x$ for some
scalar $ \xi_{x,f,\lambda}$. This implies that $x$ is an eigenvector
for every $T\in$Alg$\mathcal N$, which is imposable.  Hence we must
have $\beta_{x,f}=0$ for all $x\otimes f\in {\rm Alg}{\mathcal N}.$
Then by Eq.(3.5), we have
$$\delta(x\otimes f)=x\otimes Cf+Bx\otimes f\quad{\rm holds\  for \ all}\quad x\otimes
f\in {\rm Alg}{\mathcal N}.$$ Now, for any $T,x\otimes
f\in$Alg$\mathcal N$, by Eq.(3.9), we get
$$\begin{array}{rl}
0=&(1-\alpha)^{-1}\delta(\lambda I-T)(\lambda I-T)^{-1}x\otimes
(\lambda I-T^*)^{-1}f\\
&-(1-\alpha)^{-1}Bx\otimes (\lambda I-T^*)^{-1}f-x\otimes(\lambda
I-T^*)^{-1}Cf\\
&+x\otimes f(\lambda I-T)^{-1}\delta(\lambda I-T)(\lambda
I-T)^{-1}+x\otimes C(\lambda I-T^*)^{-1}f\\
&+(1-\alpha)^{-1}(\lambda I-T)B(\lambda I-T)^{-1}x\otimes (\lambda
I-T^*)^{-1}f,
\end{array}$$
and hence
$$\begin{array}{rl}0=&(1-\alpha)^{-1}(-\delta(T)(\lambda
I-T)^{-1}-B+(\lambda I-T)B(\lambda I-T)^{-1})x\otimes (\lambda I-T^{*})^{-1}f\\
&+x\otimes((\lambda I-T^{*})^{-1}C-(\lambda
I-T^{*})^{-1}\delta(T)^{*}(\lambda I-T^{*})^{-1}+C(\lambda
I-T^{*})^{-1})f.\end{array} $$
 This implies that $[-\delta(T)(\lambda
I-T)^{-1}-B+(\lambda I-T)B(\lambda I-T)^{-1}]x\in {\rm span} \{x\}$
for each $x\in\mathcal D_1(\mathcal N)$. So there is a scalar
$p_{\lambda}$ such that
$$\delta(T)(\lambda I-T)^{-1}+B-(\lambda I-T)B(\lambda
I-T)^{-1}=p_{\lambda} I $$ on $\mathcal D_1(\mathcal N)$. It follows
 that
$\delta(T)=BT-TB+p_{\lambda}(\lambda I-T)$ on $\mathcal D_1(\mathcal
N)$. By taking different $\lambda$, we see that $p_{\lambda}=0$, and
consequently,
$$\delta(T)|_{\mathcal D_1(\mathcal N)}=BT|_{\mathcal D_1(\mathcal N)}-TB\quad{\rm holds\  for\  all}\quad T\in {\rm Alg}{\mathcal N}.\eqno(3.14)$$

Now for any $T,S\in {\rm Alg}{\mathcal N}$, by Eq.(3.14), we have
$$\delta(TS)|_{\mathcal D_1(\mathcal N)}
=BTS|_{\mathcal D_1(\mathcal N)}-TSB =BT|_{\mathcal D_1(\mathcal
N)}S|_{\mathcal D_1(\mathcal N)}-TSB$$ and
$$\begin{array}{rl}
(\delta(T)S+T\delta(S))|_{\mathcal D_1(\mathcal N)} =&(BT|_{\mathcal
D_1(\mathcal N)}-TB)S|_{\mathcal D_1(\mathcal N)}
+T(BS|_{\mathcal D_1(\mathcal N)}-SB)\\
=&BT|_{\mathcal D_1(\mathcal N)}S|_{\mathcal D_1(\mathcal N)}-TSB.
\end{array}$$
Comparing the above two equations gives $\delta(T S)|_{\mathcal
D_1(\mathcal N)}=(\delta(T)S+T\delta(S))|_{\mathcal D_1(\mathcal
N)}$ holds for all $T,S\in {\rm Alg}{\mathcal N}$. Thus
$\delta(TS)=\delta(T)S+T\delta(S)$ holds for all $T,S\in {\rm
Alg}{\mathcal N}$ since ${\mathcal{D}_1}$ is dense in $X$.
Therefore, $\delta$ is a derivation.

{\bf Claim 3.2.}   $\delta$ is a derivation if $Z$ is an injective
operator.

In this case, by Lemma 2.5, Eqs.(3.2), (3.4)  and the equation
$$Z=[Z(I+(1-\alpha )^{-1}(\lambda
I-T)^{-1}x\otimes f)(\lambda I-T)^{-1}](\lambda I-T-x\otimes f),$$
using a similar argument to that of Claim 3.1, one can show that
$\delta$ is a derivation.

The proof of Theorem 1.1 is completed. \hfill$\Box$


\end{document}